\documentclass[12pt,a4paper,reqno]{amsart}
\usepackage{amsmath}
\usepackage{amsfonts}
\usepackage{amssymb}

\numberwithin{equation}{section}

\theoremstyle{plain}

\newtheorem{theoremenv}[subsection]{Theorem}
\newtheorem{propositionenv}[subsection]{Proposition}
\newtheorem{lemmaenv}[subsection]{Lemma}

\theoremstyle{definition}

\renewcommand{\leq}{\leqslant}
\renewcommand{\geq}{\geqslant}

\begin{document}

\newsavebox{\proofbox}
\savebox{\proofbox}{\begin{picture}(7,7)%
  \put(0,0){\framebox(7,7){}}\end{picture}}

\def\proof{\noindent\textit{Proof. }}
\def\endproof{\hfill{\usebox{\proofbox}}}
\newcommand{\Supp}{\mbox{Supp}}
\def\boxeq{\tag*{\usebox{\proofbox}}}

\def\mmm{\lambda_{b,m,N}^{(Q)}}
\def\mm{\lambda_{b,m,N}}
\def\mmmhat{\lambda_{b,m,N}^{(Q)\wedge}}
\def\mmhat{\lambda_{b,m,N}^{\wedge}}
\newcommand{\md}[1]{\ensuremath{(\mbox{mod}\, #1)}}
\newcommand{\mdsub}[1]{\ensuremath{(\mbox{\scriptsize mod}\, #1)}}
\newcommand{\mdlem}[1]{\ensuremath{(\mbox{\emph{mod}}\, #1)}}
\newcommand{\mdsublem}[1]{\ensuremath{(\mbox{\scriptsize \emph{mod}}\, #1)}}

\title{Roth's Theorem in the Primes
}

\author{Ben Green}
\address{
}
\email{bjg23@hermes.cam.ac.uk}

\thanks{The author is supported by a Fellowship of Trinity College, and for some of the period during which this work was carried out enjoyed the hospitality of Microsoft Research, Redmond WA and the Alfr\'ed R\'enyi Institute of the Hungarian Academy of Sciences, Budapest. He was supported by the \textit{Mathematics in Information Society}
project carried out by R\'enyi Institute, in the framework of the European Community's \textit{Confirming the International R\^ole of Community Research} programme.}

\begin{abstract}
We show that any set containing a positive proportion of the primes contains a 3-term arithmetic progression. An important ingredient is a proof that the primes enjoy the so-called Hardy-Littlewood majorant property. We derive this by giving a new proof of a rather more general result of Bourgain which, because of a close analogy with a classical argument of Tomas and Stein from Euclidean harmonic analysis, might be called a \emph{restriction theorem for the primes.}
\end{abstract}

\maketitle

\section{Introduction} Arguably the second most famous result of Klaus Roth is his 1953 upper bound \cite{Roth} on $r_3(N)$, defined 17 years previously by Erd\H{o}s and Tur\'an to be the density of the largest set $A \subseteq [N]$ containing no non-trivial 3-term arithmetic progression (3AP). Roth was the first person to show that $r_3(N) = o(1)$. In fact, he proved the following quantitative version of this statement.
\begin{propositionenv}[Roth] We have $r_3(N) \ll 1/\log \log N$.\endproof\end{propositionenv}
There was no improvement on this bound for nearly 40 years, until Heath-Brown \cite{HeathBrown} and Szemer\'edi \cite{Szemeredi} proved that $r_3 \ll (\log N)^{-c}$ for some small positive constant $c$. Recently Bourgain \cite{Bourgain} provided the best bound currently known.
\begin{propositionenv}[Bourgain]\label{BT} We have $r_3(N) \ll \left(\log\log N/\log N\right)^{1/2}$.\endproof\end{propositionenv}
The methods of Heath-Brown, Szemer\'edi and Bourgain may be regarded as (highly non-trivial) refinements of Roth's technique. There is a feeling that Proposition \ref{BT} is close to the natural limit of this method. This is irritating, because the sequence of primes is not covered by these results. However it \textit{is} known that the primes contain infinitely many 3APs.\footnote{In April 2004 the author and T. Tao published a preprint shoing that the primes contain arbitrarily long arithmetic progressions.}
\begin{propositionenv}[Van der Corput] The primes contain infinitely many \emph{3AP}s.\endproof
\end{propositionenv}
Van der Corput's method is very similar to that used by Vinogradov to show that every large odd number is the sum of three primes. Let us also mention a paper of Balog \cite{Balog} in which it is shown that for any $n$ there are $n$ primes $p_1,\dots,p_n$ such that \textit{all} of the averages $\frac{1}{2}(p_i + p_j)$ are prime.
In this paper we propose to prove a common generalization of the results of Roth and Van der Corput. Write $\mathcal{P}$ for the set of primes.
\begin{theoremenv}\label{mainthm} Every subset of $\mathcal{P}$ of positive upper density contains a \emph{3AP}. 
\end{theoremenv}
In fact, we get an explicit upper bound on the density of a 3AP-free subset of the primes, but it is ridiculously weak. Observe that as an immediate consequence of Theorem \ref{mainthm} we obtain what might be termed a van der Waerden theorem in the primes, at least for progressions of length 3. That is, if one colours the primes using finitely many colours then one may find a monochromatic 3AP.

We have not found a written reference for the question answered by Theorem \ref{mainthm}, but M.N. Huxley has discussed it with several people \cite{Huxley}.

To prove Theorem \ref{mainthm} we will use a variant of the following result. This says that the primes enjoy what is known as the Hardy-Littlewood majorant property. 
\begin{theoremenv}\label{thmHLM}
Suppose that $p \geq 2$ is a real number, and let $\mathcal{P}_N = \mathcal{P} \cap [1,N]$. Let $\{a_n\}_{n \in \mathcal{P}_N}$ be any sequence of complex numbers with $|a_n| \leq 1$ for all $n$. Then
\begin{equation}\label{HLM} \left\Vert \sum_{n \in \mathcal{P}_N} a_n e(n\theta) \right\Vert_{L^p(\mathbb{T})}  \leq  C(p)\left\Vert \sum_{n \in \mathcal{P}_N} e(n\theta) \right\Vert_{L^{p}(\mathbb{T})},\end{equation}
where the constant $C(p)$ depends only on $p$.
\end{theoremenv}
It is perhaps surprising to learn that such a property does not hold with \textit{any} set $\Lambda \subseteq [N]$ in place of $\mathcal{P}_N$. Indeed, when $p$ is an even integer it is rather straightforward to check that any set \textit{does} satisfy \eqref{HLM} (with $C(p) = 1$). However, there are sets for which \eqref{HLM} fails badly when $p$ is not an even integer. For a discussion of this see \cite{GreenRuzsaHLM} and for related matters including connections with the Kakeya problem, see \cite{MockenhauptSchlag,Mockenhaupt}.

We will apply a variant of Theorem \ref{thmHLM} for $p = 5/2$, when it certainly does not seem to be trivial. To prove it, we will establish a somewhat stronger result which we call a \textit{restriction theorem for primes}. The reason for this is that our argument is very closely analogous to an argument of Tomas and Stein \cite{Tomas} concerning Fourier transforms of measures supported on spheres.

A proof of the restriction theorem for primes was described, in a different context, by Bourgain \cite{B1}. Our argument, being visibly analagous to the approach of Tomas, is different and has more in common with \S 3 of \cite{B2}. This more recent paper of Bourgain deals with restriction phenomena of certain sets of lattice points.

To deduce Theorem \ref{mainthm} from (a variant of) Theorem \ref{thmHLM} we use a variant of the technique of \textit{granularization} as developed by I.Z. Ruzsa and the author in a series of papers beginning with \cite{GreenRuzsa}, as well as a ``statistical'' version of Roth's theorem due to Varnavides. We will also require an argument of Marcinkiewicz and Zygmund which allows us to pass from the continuous setting in results such as \eqref{HLM} -- that is to say, $\mathbb{T}$ -- to the discrete, namely $\mathbb{Z}/N\mathbb{Z}$.

Finally, we would like to remark that it is possible, indeed probable, that Roth's theorem in the primes is true on grounds of density alone. The best known lower bound on $r_3(N)$ comes from a result of Behrend \cite{Behrend} from 1946.
\begin{propositionenv}[Behrend] We have $r_3(N) \geq e^{-C\sqrt{\log N}}$ for some absolute constant $C$.\endproof\end{propositionenv}
This may well give the correct order of magnitude for $r_3(N)$, and if anything like this could be proved Theorem \ref{mainthm} would of course follow trivially.
\section{Preliminaries and an outline of the argument} Although the main results of this paper concern the primes in $[N]$, it turns out to be necessary to consider slightly more general sets. Let $m \leq \log N$ be a positive integer and let $b$, $0 \leq b \leq m - 1$, be coprime to $m$. We may then define a set
\[ \Lambda_{b,m,N}  =  \left\{n \leq N\, |\, nm + b \; \mbox{is prime}\right\}.\]
We expect $\Lambda_{b,m,N}$ to have size about $mN/\phi(m)\log N$, and so it is natural to define a function $\mm$ supported on $\Lambda_{b,m,N}$ by setting
\[ \lambda_{b,m,N}(n)  =  \left\{\begin{array}{ll} \phi(m)\log (nm + b)/mN & \mbox{if $n \in \Lambda_{b,m,N}$}\vspace{2mm} \\ 0 & \mbox{otherwise}.\end{array}\right.\]
For simplicity we write $X = \Lambda_{b,m,N}$ for the next few pages. We will abuse notation and consider $\lambda_{b,m,N}$ as a measure on $X$. Thus for example $\mm(X)$, which is defined to be $\sum_n \mm(n)$, is roughly $1$ by the prime number theorem in arithmetic progressions. We use $L^p(d\mm)$ norms and also the inner product $\langle f,g \rangle_X = \sum f(n)\overline{g(n)}\mm(n)$ without further comment. 

It is convenient to use the wedge symbol for the Fourier transforms on both $\mathbb{T}$ and $\mathbb{Z}$, which we define by $f^{\wedge}(n)  = \int f(\theta)e(-n\theta)\,d\theta$ and $g^{\wedge}(\theta) = \sum_n g(n)e(n\theta)$ respectively. Here, of course, $e(\alpha) = e^{2\pi i \alpha}$.

For any measure space $Y$ let $B(Y)$ denote the space of continuous functions on $Y$ and define a map $T : B(X) \rightarrow B(\mathbb{T})$ via
\begin{equation}\label{eq608}T : f \longmapsto (f\lambda_{b,m,N})^{\wedge}.\end{equation} 
The object of this section is to give a new proof of the following result, which may be a called a \textit{restriction theorem for primes}. 
\begin{theoremenv}[Bourgain]\label{primerest} Suppose that $p > 2$ is a real number. Then there is a constant $C(p)$ such that for all functions $f : X \rightarrow \mathbb{C}$ we have
\begin{equation}\label{eq1} \Vert Tf \Vert_{p}  \leq  C(p)N^{-1/p}\Vert f \Vert_{2}.\end{equation}
\end{theoremenv}
Remember that the $L^2$ norm is taken with respect to the measure $\mm$.
Theorem \ref{primerest} probably has most appeal when $b = m = 1$, in which case we may derive consequences for the primes themselves. Later on, however, we will take $m$ to be a product of small primes, and so it is necessary to have the more general form of the theorem.

We turn now to an outline of the proof of Theorem \ref{primerest}. The analogy between our proof and an argument by Tomas \cite{Tomas}, giving results of a similar nature for spheres in high-dimensional Euclidean spaces, is rather striking. In fact, the reader may care to look at the presentation of Tomas's proof in \cite{Tao}, whereupon she will see that there is an almost exact correspondence between the two arguments.

To begin with, the proof proceeds by the method of $T$ and $T^{\ast}$, a basic technique in functional analysis. One can check that the operator $T^{\ast} : B(\mathbb{T}) \rightarrow B(X)$ is given by
\begin{equation}\label{eq027} T^{\ast} : g \longmapsto g^{\wedge}|_{X},\end{equation} by verifying the relation
\[ \langle Tf,g\rangle_{\mathbb{T}}  =  \int (f\mm)^{\wedge}(\theta)\overline{g(\theta)}\,d\theta  =  \sum_n f(n)\overline{g^{\wedge}(n)}\mm(n)  =  \langle f, T^{\ast}g\rangle_X.\]
The equation \eqref{eq027} explains the term \textit{restriction}.
Using \eqref{eq027} we see that the operator $TT^{\ast}$ is the map from $B(\mathbb{T})$ to itself given by
\begin{equation}\label{TastT} TT^{\ast} : f \longmapsto f \ast \mm^{\wedge}.\end{equation}
Now Theorem \ref{primerest} may be written, in an obvious notation, as 
\begin{equation}\label{p96} \Vert T \Vert_{2 \rightarrow p}  \leq  C(p)N^{-1/p}.\end{equation}
The principle of $T$ and $T^{\ast}$, as we will use it, states that
\begin{equation}\label{TTstar} \Vert T \Vert_{2 \rightarrow p}^2 
   =  \Vert TT^{\ast} \Vert_{p' \rightarrow p}  =  \Vert T^{\ast} \Vert_{p' \rightarrow 2}^2 \end{equation}
We would like to emphasise that there is nothing mysterious going on here -- this result is just an elegant and convenient way of bundling together some applications of H\"older's inequality. The proof of the part that we will need, that is to say is the inequality $\Vert T \Vert_{2 \rightarrow p}^2 
 \leq \Vert TT^{\ast} \Vert_{p' \rightarrow p}$, is simply
\begin{eqnarray*} \Vert Tf \Vert_p & = & \sup_{\Vert g \Vert_{p'} = 1} \langle Tf,g\rangle\\ & = & \sup_{\Vert g \Vert_{p'} = 1} \langle f,T^{\ast}g\rangle \\ & \leq & \Vert f\Vert_2 \sup_{\Vert g\Vert_{p'} = 1}\Vert T^{\ast} g\Vert_2 \\ & = &  \Vert f\Vert_2 \sup_{\Vert g\Vert_{p'} = 1}\langle g,TT^{\ast}g\rangle^{1/2} \\ & \leq & \Vert f \Vert_2\Vert TT^{\ast}\Vert_{p' \rightarrow p}^{1/2}. \end{eqnarray*}
Thus we will, for much of the paper, be concerned with showing that the operator $TT^{\ast}$ as given by \eqref{TastT} satisfies the bound
\begin{equation}\label{satis} \Vert TT^{\ast} \Vert_{p' \rightarrow p}  \leq  C'(p)N^{-2/p}.\end{equation}
The preceding remarks show that a proof of this will imply Theorem \ref{primerest}.
To get such a bound one splits $\lambda$ into certain \textit{dyadic pieces}, that is a sum \begin{equation}\label{dyadic} \lambda_{b,m,N}  =  \sum_{j = 1}^K \psi_j + \psi_{K+1}.\end{equation} The slightly curious way of writing this indicates that the definition of $\psi_{K+1}$ will be a little different from that of the other $\psi_j$. We will define these pieces so that they satisfy the $L^1$--$L^{\infty}$ estimates 
\begin{equation}\label{eq699} \Vert f \ast \psi_j^{\wedge} \Vert_{\infty}  \ll_{\epsilon}  2^{-(1 - \epsilon)j} \Vert f \Vert_{1}\end{equation} for some $\epsilon < (p-2)/2$, and also the $L^2$--$L^2$ estimates
\begin{equation}\label{eq698} \Vert f \ast \psi_j^{\wedge} \Vert_{2}  \ll_{\epsilon}  \frac{2^{\epsilon j}}{N} \Vert f \Vert_2.\end{equation}
Applying the Riesz-Thorin interpolation theorem (see \cite{GreenRKP}, Chapter 7) will then give
\[ \Vert f \ast \psi_j^{\wedge} \Vert_{p}  \ll  2^{-\delta j}N^{-2/p} \Vert f \Vert_{p'}\] for some positive $\delta$ (depending on $\epsilon$). Summing these estimates from $j = 1$ to $K+1$ will establish \eqref{satis} and hence Theorem \ref{primerest}.

To define the decomposition \eqref{dyadic} we need yet more notation. From the outset we will suppose that we are trying to prove Theorem \ref{primerest} for a particular value of $p$ -- the argument is highly and essentially non-uniform in $p$. Write $A = 4/(p-2)$.  Let $1 < Q \leq (\log N)^A$. If $b,m,N$ are as before (recall that $m \leq \log N$) then we define a measure $\mmm$ on $\mathbb{Z}$ by setting
\[ \mmm(n)  =  \left\{\begin{array}{ll} N^{-1}\prod_{\substack{p \leq Q \\p\nmid m}} \left(1 - \frac{1}{p}\right)^{-1} & \mbox{if $n\leq N$ and $p \,|\,(nm + b) \Rightarrow p > Q$}\vspace{2mm} \\
 0 & \mbox{otherwise}.\end{array}\right.\]
Define $\lambda_{b,m,N}^{(1)}(n) = 0$ for all $n$.

As $Q$ becomes large the measures $\mmm$ look more and more like $\mm$. Much of \S \ref{sec4} will be devoted to making this principle precise. We will sometimes refer to the support of $\mmm$ as the set of $Q$-\textit{rough} numbers.

Now let $K$ be the smallest integer with
\begin{equation}\label{Kdef} 2^K  >  \textstyle\frac{1}{10}\displaystyle(\log N)^{A}\end{equation}
and define 
\begin{equation}\label{eq681} \psi_j  =  \mm^{(2^{j})} - \mm^{(2^{j-1})}\end{equation} for\ $j = 1,\dots,K$ and define \begin{equation}\label{eq6822}\psi_{K+1}  =  \mm - \mm^{(2^{K})},\end{equation} so that \eqref{dyadic} holds. In the next two sections we prove the two required estimates, \eqref{eq699} and \eqref{eq698}. 

Let us note here that the main novelty in our proof of Theorem \ref{primerest} lies in the definition of the dyadic decomposition \eqref{dyadic}. By contrast, the analogous dyadic decompositions in \cite{B2} take place on the Fourier side, requiring the introduction of various smooth cutoff functions not specifically related to the underlying arithmetic structure.
\section{An $L^2$--$L^2$ estimate} It turns out that the proof of \eqref{eq698}, the $L^2$--$L^2$ estimate, is by far the easier of the two estimates required. We have
\begin{eqnarray*} \Vert f \ast \psi_j^{\wedge}\Vert_2 & = & \Vert \widehat{f} \psi_j \Vert_2 \\ & \leq & \Vert \psi_j \Vert_{\infty} \Vert \widehat{f} \Vert_2 \\ & = & \Vert \psi_j \Vert_{\infty} \Vert  f \Vert_2.\end{eqnarray*}
Suppose first of all that $1 \leq j \leq K$. Then
\begin{eqnarray*} \Vert\psi_j\Vert_{\infty} & \leq & \Vert \mm^{(2^{j})}\Vert_{\infty} + \Vert \mm^{(2^{j-1})}\Vert_{\infty} \\
& = & N^{-1}\prod_{\substack{p \leq 2^{j+1} \\ p\nmid m}} \left(1 - \frac{1}{p}\right)^{-1} + N^{-1}\prod_{\substack{p \leq 2^{j} \\ p\nmid m}} \left(1 - \frac{1}{p}\right)^{-1}.\end{eqnarray*}
The two products here may be estimated using Merten's formula \cite[Chapter 22]{HardyWright}: \[ \prod_{p \leq Q}(1 - p^{-1})  \sim  \frac{e^{-\gamma}}{\log Q}.\] This gives
\begin{equation}\Vert \psi_j \Vert_{\infty}  \ll  j/N,\end{equation} and hence
\begin{equation}\label{later1} \Vert f \ast \psi_j^{\wedge}\Vert_2  \ll  \frac{j}{N}\Vert f \Vert_2,\end{equation}
which is certainly of the requisite form \eqref{eq698}. For $j = K+1$ we have
\begin{eqnarray*}
\Vert \psi_{K+1} \Vert_{\infty} & \leq & \Vert \mm^{(2^K)} \Vert_{\infty} + \Vert \mm \Vert_{\infty} \\ & \ll & \log N/N,\end{eqnarray*}
so that
\begin{equation}
\label{later22}
\Vert f \ast \psi_{K+1}^{\wedge}\Vert_2  \ll  \frac{\log N}{N}\Vert f \Vert_2. \end{equation}
This also constitutes an estimate of the type \eqref{eq698} for some $\epsilon < (p-2)/2$. Indeed, recalling our choice of $A$ and $K$ (viz. \eqref{Kdef}) one can check that $2^K \geq (\log N)^{1/\epsilon}$ for some such $\epsilon$.
\section{An $L^1$--$L^{\infty}$ estimate}\label{sec4} This section is devoted to the rather lengthy task of proving estimates of the form \eqref{eq699}.\vspace{11pt}

\noindent\textit{Introduction.}
The first step towards obtaining an estimate of the form \eqref{eq699} is to observe that
\begin{equation}\label{later2} \Vert f \ast \psi_j^{\wedge} \Vert_{\infty}  \leq  \Vert \psi_j^{\wedge} \Vert_{\infty}\Vert f \Vert_1.\end{equation}
We will prove that $\Vert \psi_j^{\wedge} \Vert_{\infty}$ is not too large by proving
\begin{propositionenv}\label{mainprop} Suppose that $Q \leq (\log N)^A$. Then we have the estimate \[\Vert \mm^{\wedge} - \mm^{(Q)\wedge}\Vert_{\infty}  \ll  \log \log Q/Q.\]
\end{propositionenv}
The detailed proof of this fact will occupy us for several pages. Let us begin, however, by using \eqref{later2} to see how it implies an estimate of the form \eqref{eq699}. If $1 \leq j \leq K$ then we have
\begin{eqnarray}\nonumber
\Vert \psi_j^{\wedge} \Vert_{\infty} & = & \Vert \mm^{(2^{j})\wedge} - \mm^{(2^{j-1})\wedge}\Vert_{\infty} \\ \nonumber & \leq & \Vert \mm^{\wedge} - \mm^{(2^{j})\wedge}\Vert_{\infty} + \Vert \mm^{\wedge} - \mm^{(2^{j-1})\wedge}\Vert_{\infty} \\ \label{later3}& \leq & \log j/2^j.
\end{eqnarray}
This is certainly of the form \eqref{eq699}. The estimate for $j = K+1$ is even easier, being immediate from Proposition \ref{mainprop}.

To prove Proposition \ref{mainprop} we will use the Hardy-Littlewood circle method.
Thus we divide $\mathbb{T}$ into two sets, traditionally referred to as the major and minor arcs. It is perhaps best if we define these explicitly at the outset. Thus let $p$ be the exponent for which we are trying to prove Theorem \ref{primerest}. Recall that $A = 4/(p-2)$, and set $B = 2A + 20$. These numbers will be fixed throughout the proof. By Dirichlet's theorem on approximation, every $\theta \in \mathbb{T}$ satisfies 
\begin{equation}\label{u72} \left|\theta - \frac{a}{q}\right| \leq  \frac{(\log N)^B}{qN}\end{equation} for some $q \leq N(\log N)^{-B}$ and some $a$, $(a,q) = 1$. 
The major arcs consist of those $\theta$ for which $q$ can be taken to be at most $(\log N)^B$. We will write this collection using the notation
\[ \mathfrak{M}  =  \bigcup_{\substack{q \leq (\log N)^B \\ (a,q) = 1}} \mathfrak{M}_{a,q}.\] For these $\theta$, the Fourier transforms $\mm^{(Q)\wedge}$ and $\mm^{\wedge}$ depend on the distribution of the almost-primes and primes along arithmetic progressions with common difference at most $(\log N)^B$. The minor arcs $\mathfrak{m}$ consist of all other $\theta$. Here different techniques apply, and one can conclude that both $\mm^{(Q)\wedge}$ and $\mm^{\wedge}$ are small. The triangle inequality then applies.

The ingredients are as follows. The almost-primes are eminently suited to applications of sieve techniques. To keep the paper as self-contained as possible, we will follow Gowers \cite{Gowers} and use arguably the simplest sieve, that due to Brun, on both the major and minor arcs.

 The genuine primes, on the other hand, are harder to deal with. Here we will quote two well-known results from the literature. The information concerning distribution along arithmetic progressions to small moduli comes from the prime number theorem of Siegel and Walfisz. 
\begin{propositionenv}[Siegel--Walfisz] \label{Sieg} Suppose that $q \leq (\log N)^B$, that $(a,q) = 1$ and that $1 \leq N_1 \leq N_2 \leq N$. Then
\begin{equation}\boxeq \sum_{\substack{N_1 < p \leq N_2 \\ p \equiv a \mdsublem{q}}} \log p  =  \frac{N_2 - N_1}{\phi(q)} + O\left(N\exp(-C_B\sqrt{\log N})\right).\end{equation}
\end{propositionenv}
The rather strange formulation of the theorem reflects the fact that the constant $C_B$ is \textit{ineffective} for any $B \geq 1$ due to the possible existence of a \textit{Siegel zero}. For more information, including a complete proof of Proposition \ref{Sieg}, see Davenport's book \cite{davenport}.

The techniques for dealing with the minor arcs are associated with the names of Weyl, Vinogradov and Vaughan.\vspace{11pt}

\noindent\textit{The major arcs.}
We will have various functons $f : [N] \rightarrow \mathbb{R}$ with
\begin{equation}\label{assump2} \Vert f \Vert_{\infty}  =  O(\log N/N)\end{equation} which are regularly distributed along arithmetic progressions in the following sense. If $L \geq N(\log N)^{-2B - A - 1}$ and if $X \subseteq [N]$ is an arithmetic progression $\{r,r+q,\dots,r + (L-1)q\}$ with $q \leq (\log N)^B$ then
\begin{equation}\label{assump} \sum_{n \in X} f(n)  =  \frac{L}{N}\left(\gamma_{r,q}(f) + O((\log N)^{-A})\right),\end{equation}
where $\gamma_{r,q}$ depends only on $r$ and $q$, $|\gamma_{r,q}| \leq q$ and the implied constant in the $O$ term is absolute. This information is enough to get asymptotics for $f^{\wedge}(\theta)$ when $|\theta - a/q|$ is small, as we prove in the next few lemmas.

For a residue $r$ modulo $q$, write $N_r$ for the set $\{n \leq N : n \equiv r \md{q}\}$. Write $\tau$ for the function on $\mathbb{T}$ defined by $\tau(\theta) = N^{-1}\sum_{n \leq N} e(\theta n)$. The first lemma deals with $f^{\wedge}(\theta)$ for $|\theta| \leq (\log N)^B/qN$.
\begin{lemmaenv}\label{lem5}
Let $r$ be a residue modulo $q$, suppose that $|\theta| \leq (\log N)^B/qN$, and suppose that the function $f$ satisfies \eqref{assump2} and \eqref{assump}. Then 
\[ \sum_{n \in N_r} f(n)e(\theta n)  =  q^{-1}\gamma_{r,q}(f)\tau(\theta) + O(q^{-1}(\log N)^{-A}).\] 
\end{lemmaenv}
\proof Set $L = N(\log N)^{-2B - A - 1}$ and partition $N_r$ into arithmetic progressions $(X_i)_{i = 1}^T$ of common difference $q$ and length between $L$ and $2L$, where $T \leq 2N/Lq$. For each $i$ fix an element $x_i \in X_i$. 
\begin{eqnarray}\nonumber
\sum_{n \in N_r} f(n)e(\theta n) & = & \sum_{i = 1}^T \sum_{n \in X_i} f(n)e(\theta n) \\ \nonumber & = & \sum_{i = 1}^T e(\theta x_i)\sum_{n \in X_i} f(n) + \sum_{i = 1}^T \sum_{n \in X_i}f(n)\left(e(\theta n) - e(\theta x_i)\right) \\ \nonumber & = & \sum_{i = 1}^T e(\theta x_i)\frac{|X_i|}{N}\left(\gamma_{r,q}(f) + O((\log N)^{-A})\right) + O(LN^{-1}q^{-1}(\log N)^{B+1}) \\ \label{e723}& = & \gamma_{r,q}(f)\sum_{i=1}^T e(\theta x_i)\frac{|X_i|}{N} + O\left(q^{-1}(\log N)^{-A}\right).\end{eqnarray} However
\begin{eqnarray}\nonumber
\sum_{i=1}^T e(\theta x_i)|X_i| & = & \sum_{i = 1}^T \sum_{n \in X_i} e(\theta n) + \sum_{i = 1}^T \sum_{n \in X_i}\left(e(\theta x_i) - e(\theta n)\right) \\ \label{e724}& = & \sum_{n \in N_r} e(n\theta) + O(Lq^{-1}(\log N)^B).
\end{eqnarray}
Finally, observe that if $0 \leq r,s \leq q-1$ then
\[ \sum_{n \in N_r} e(\theta n) - \sum_{n \in N_s} e(\theta n)  =  O((\log N)^B),\] and so
\[
\left|N^{-1}\sum_{n \in N_{r}} e(\theta n) - q^{-1}\tau(\theta)\right|  =  O(N^{-1}(\log N)^B).\] Combining this with \eqref{e723} and \eqref{e724} completes the proof of the lemma.\endproof\vspace{11pt}

We may now get an asymptotic for $f^{\wedge}(\theta)$ when $\theta$ is in the neighbourhood of $a/q$.
\begin{lemmaenv}\label{tizenketto}
Suppose that $f$ satisfies the conditions \eqref{assump2} and \eqref{assump} and that $\theta \in \mathfrak{M}_{a,q}$ for some $a,q$ with $(a,q) = 1$ and $q \leq (\log N)^B$. Write
\begin{equation}\label{61.1} \sigma_{a,q}(f)  =  \sum_r e(ar/q)\gamma_{r,q}(f).\end{equation}
Then we have
\begin{equation}\label{61.2} f^{\wedge}(\theta)  =  q^{-1}\sigma_{a,q}(f)\tau(\theta - a/q) + O((\log N)^{-A}).\end{equation}
\end{lemmaenv}
\proof Write $\beta = \theta - a/q$. Then
\begin{eqnarray*}
f^{\wedge}(\theta) & = & \sum_{n \leq N}f(n)e(\theta n) \\ & = & \sum_{r\mdsub{q}} e(ar/q) \sum_{n \in N_r}f(n) e(\beta n) \\ & = & q^{-1}\tau(\beta)\sum_{r\mdsub{q}} e(ar/q)\gamma_{r,q}(f) + O((\log N)^{-A})\\ & = & q^{-1}\sigma_{a,q}(f)\tau(\beta) + O((\log N)^{-A}).\end{eqnarray*}
This concludes the proof of the lemma.\endproof\vspace{11pt}

To apply these lemmas, we need to show that $f = \mm^{(Q)}$ and $f = \mm$ satisfy \eqref{assump2} and \eqref{assump} for suitable choices of $\gamma_{r,q}(f)$. We will then evaluate the sums $\sigma_{a,q}(f)$. This slightly tedious business is the subject of our next four lemmas.
\begin{lemmaenv}\label{lem23} $f = \mm$ satisfies \eqref{assump2} and \eqref{assump} with 
\[ \gamma_{r,q}(f)  =  \left\{\begin{array}{ll} \phi(m)q/\phi(mq) & \mbox{if $(mr + b,mq) = 1$}\vspace{2mm} \\ 0 & \mbox{otherwise}.\end{array}\right.\]
\end{lemmaenv}
\proof This is a fairly immediate consequence of the Siegel-Walfisz Theorem (Proposition \ref{Sieg}). Let $X = \{r,r+q,\dots,r + (L-1)q\}$ be any progression contained in $[N]$ with common difference $q \leq (\log N)^B$ and length $L \geq N(\log N)^{-2B - A - 1}$. An element $r + jq \in X$ lies in $\Lambda_{b,m,N}$ precisely if $(mr + b) + jmq$ is prime, so the lemma is trivially true unless $(mr + b,mq) = 1$. Supposing this to be the case, we may use Proposition \ref{Sieg}. Recalling that $m \leq \log N$, one has
\begin{eqnarray*}
\mm(X) & = & \frac{\phi(m)qL}{\phi(mq)N} +  O\left(mq\exp(-C_{B+1} \sqrt{\log mqN})\right)\\ & = & 
\frac{L}{N}\left(\frac{\phi(m)q}{\phi(mq)} + O((\log N)^{-A})\right),\end{eqnarray*}
as required.\endproof
\begin{lemmaenv}\label{lem24} $f = \mmm$ satisfies \eqref{assump2} and \eqref{assump} with
\[ \gamma_{r,q}(f)  =  \left\{\begin{array}{ll} \prod_{\substack{p \leq Q \\ p \nmid m}}\left(1 - \frac{1}{p}\right)^{-1}\prod_{\substack{p \leq Q \\p \nmid mq}} \left(1 - \frac{1}{p}\right) & \mbox{if $(mr + b,mq)$ is $Q$-rough}\vspace{2mm} \\ 0 & \mbox{otherwise}.\end{array}\right.\]
\end{lemmaenv}
\proof Consider an arithmetic progression $X = \{r,r+q,\dots,r+(L-1)q\}$. Let $p_1,\dots,p_k$ be the primes with $p \leq Q$ and $p \nmid m$. If $(mr +b,mq)$ is not $Q$-rough then $p_i | (mr + b,mq)$ for some $i$, and the second alternative of the lemma clearly holds. Suppose then that $(mr + b,mq)$ is $Q$-rough. We will apply the Brun sieve to estimate $\mmm(X)$.

Let $x \in X$ be chosen uniformly at random, and for each $i$ let $X_i$ be the event $p_i | (mx + b)$. Since $p_i \nmid(mr + b,mq)$, the probability of $X_i$ is $\epsilon_i/p_i + O(L^{-1})$, where $\epsilon_i = 0$ if $p_i | q$ and $\epsilon_i = 1$ otherwise.
Now we have
\begin{equation}\label{gg290} \frac{N}{L}\prod_{\substack{p \leq Q \\p \nmid m}}\left(1 - \frac{1}{p}\right)\mmm(X)  =  \mathbb{P}\left(\bigcap X_i^{c}\right)  =  U,\end{equation} say. By the inclusion-exclusion formula it follows that for every positive integer $t$
\begin{equation}\label{eq3} U  =  \sum_{s=0}^t (-1)^s \sum_{1 \leq i_1 < \dots < i_s \leq k} \prod_{j=1}^s \epsilon_{i_j}/p_{i_j} + O(L^{-1}) \sum_{s = 1}^t \binom{k}{s}.\end{equation}  It is helpful to have the error term here in a more usable form. To this end, observe that it is certainly at most $O(k^t/L)$. We wish to replace the main term in \eqref{eq3} by $\prod_{i=1}^k \left(1 - \epsilon_i/p_i\right)$, which is equal to the completed sum
\[ \sum_{s=0}^k (-1)^s \sum_{1 \leq i_1 < \dots < i_s \leq k} \prod_{j=1}^s \epsilon_{i_j}/p_{i_j}.\] Doing this introduces an error
\[ E  =  \sum_{s=t+1}^k (-1)^s \sum_{1 \leq i_1 < \dots < i_s \leq k} \prod_{j=1}^s \epsilon_{i_j}/p_{i_j},\] which is bounded above by
\begin{equation}\label{eq4} \sum_{s = t+1}^k \frac{1}{s!} \left(\sum_{i=1}^k \frac{1}{p_i}\right)^s.\end{equation} By another result of Mertens one has $\sum_{i=1}^k p_i^{-1} \leq \log\log Q + O(1)$. Hence if $t \geq 3\log\log Q$ then each term in \eqref{eq4} is at most one half the previous one, leading to the bound
\[ |E|  \leq  \frac{2(\log\log Q)^t}{t!}  \leq  \left(\frac{4e \log \log Q}{t}\right)^t.\] 
Combining all of this gives
\[ U  =  \prod_{i=1}^k \left(1 - \epsilon_i/p_i\right)  +  O(k^t/L)  +  O\left((4e\log\log Q/t)^t\right).\] Using the trivial bound $k \leq Q$, and choosing $t = \log N/2A\log \log N$, one gets
\begin{eqnarray*} U & = & \prod_{i=1}^k \left(1 - \epsilon_i/p_i\right) + O(N^{-1/4A})\\ & = & \prod_{\substack{p \leq Q \\ p \nmid mq}}\left(1 - \frac{1}{p}\right) + O(N^{-1/4A}).\end{eqnarray*} The lemma is immediate from this and \eqref{gg290}; we have
\begin{eqnarray*} \mmm(X) & = & \prod_{\substack{p \leq Q \\ p \nmid m}}\left(1 - \frac{1}{p}\right)^{-1} \cdot \frac{L}{N} \cdot\left( \prod_{\substack{p \leq Q \\ p \nmid mq}}\left(1 - \frac{1}{p}\right) + O(N^{-1/4A})\right) \\ & = & \frac{L}{N}\left(\gamma_{r,q} + O((\log N)^{-A})\right),\end{eqnarray*}
where $\gamma_{r,q}$ has the form claimed.\endproof\vspace{11pt}

Building on the last lemma, the next lemma gives an evaluation of $\sigma_{a,q}(\mmm)$ and an asymptotic for $\mmmhat(\theta)$ when $\theta \in \mathfrak{M}_{a,q}$. If $Q \geq 2$ we say that a positive integer is $Q$-\textit{smooth} if all of its prime divisors are at most $Q$. We declare there to be no $1$-smooth numbers.
\begin{lemmaenv}\label{lem27} Suppose that $(a,q) = 1$. Then 
\[\sigma_{a,q}(\mmm)  =  \left\{\begin{array}{ll} \displaystyle\frac{q\mu(q)}{\phi(q)}e\left(-\frac{ab\overline{m}}{q}\right) & \mbox{if $(m,q) = 1$ and $q$ is $Q$-smooth;}\vspace{2mm} \\ 0 & \mbox{otherwise},\end{array}\right.\]
where $\overline{m}$ is the inverse of $m$ modulo $q$.
If $\theta \in \mathfrak{M}_{a,q}$ then 
\[\mmmhat(\theta)  =  \left\{\begin{array}{ll} \displaystyle\frac{\mu(q)}{\phi(q)}e\left(-\frac{ab\overline{m}}{q}\right)\tau\left(\theta - \frac{a}{q}\right) + O((\log N)^{-A}) & \mbox{if $(m,q) = 1$} \\ & \qquad \mbox{and $q$ is $Q$-smooth;}\vspace{2mm} \\ O\left((\log N)^{-A}\right) & \mbox{otherwise}.\end{array}\right.\]
\end{lemmaenv}
\proof Recall the definition \eqref{61.1} of $\sigma_{a,q}$, and also Lemma \ref{lem24}. We shall prove that
\begin{equation}\label{eq444} \sum_{\substack{r \mdsub{q} \\(mr + b,mq)\, \mbox{\scriptsize is $Q$-rough}}} e(ar/q)  =  \left\{\begin{array}{ll} e(-ab\overline{m}/q)\mu(q) & \mbox{if $(m,q) = 1$ and $q$ is $Q$-smooth}\vspace{2mm} \\ 0 & \mbox{otherwise}.\end{array}\right.\end{equation}
Now if $p | m$ then $p$ can never divide $mr + b$, because we are assuming that $(m,b) = 1$. Let $q_0$ be the largest factor of $q$ which is a product of primes $p$ with $p \leq Q$ and $p \nmid m$. Then the sum \eqref{eq444} is just
\begin{equation}\label{eq472} \sum_{\substack{r \mdsub{q} \\(q_0,mr + b) = 1}} e(ar/q).\end{equation}
Set $q_1 = q/q_0$ and write, for each $r$ mod $q$, $r = kq_0 + s$ where $0 \leq k \leq q_1 - 1$ and $s$ is a residue mod $q_0$. Then the sum \eqref{eq472} is
\[
\sum_{\substack{s \mdsub{q_0} \\(q_0,ms + b) = 1}} \sum_{k = 0}^{q_1 - 1} e\left(\frac{a(kq_0 + s)}{q}\right)  =  \sum_{\substack{s \mdsub{q_0} \\(q_0,mr + b) = 1}} e(as/q)\sum_{k = 0}^{q_1 - 1} e(ak/q_1).\] Now $a$ is coprime to $q$ and hence to $q_1$, and therefore the rightmost sum here vanishes unless $q_1 = 1$. This is the case precisely if $q_0 = q$, which means that $(q,m) = 1$ and $q$ is $Q$-smooth. In this case, the sum is
\[ \sum_{\substack{s \mdsub{q} \\(q,ms + b) = 1}} e(as/q).\] Set $t = ms + b$. Then this sum is just
\begin{eqnarray*} \sum_{\substack{t \mdsub{q}\\(q,t) = 1}} e\left(\frac{a\overline{m}(t - b)}{q}\right) & = & e(-ab\overline{m}/q)\sum_{(q,t) = 1}e(a\overline{m}t/q)\\ & = & e(-ab\overline{m}/q)\mu(q).\end{eqnarray*}
This last evaluation, of what is known as a \textit{Ramanujan Sum}, is well-known and is contained, for example, in \cite{HardyWright}. This proves \eqref{eq444}.

Now to obtain $\sigma_{a,q}$ we must simply multiply \eqref{eq444} by the factor 
\[F  =   \prod_{\substack{p < Q \\ p \nmid m}}\left(1 - \frac{1}{p}\right)^{-1}\prod_{\substack{p \leq Q \\ p \nmid mq}} \left(1 - \frac{1}{p}\right) \] appearing in Lemma \ref{lem24}. One gets zero unless $(m,q) = 1$ and $q$ is $Q$-smooth, in which case it is not hard to see that $F = q/\phi(q)$. This completes the evaluation of $\sigma_{a,q}(\mmm)$, and the claimed form for $\mmmhat(\theta)$ is an immediate consequence of Lemma \ref{tizenketto}.\endproof\vspace{11pt}

We need a version of the above lemma in which $\mmm$ is replaced by $\mm$. Fortunately, we can save ourselves some work by noticing that for fixed $q$ and $m$ we have
\begin{equation}\label{eq774} \gamma_{r,q}(\mm)  =  \gamma_{r,q}(\mmm)\end{equation} for sufficiently\footnote{Here we regard $\gamma_{r,q}(\mm)$ and $\gamma_{r,q}(\mmm)$ as purely formal expressions, so there is no issue of whether or not, for example, Lemma \ref{lem27} is valid for ``sufficiently large'' $Q$.} large $Q$. Thus $\sigma_{a,q}(\mm)$ can be evaluated by simply letting $Q \rightarrow \infty$ in the first formula of Lemma \ref{lem27}. We get  
\begin{equation}\label{28.1} \sigma_{a,q}(\mm)  =  \left\{\begin{array}{ll} q\mu(q)e(-ab\overline{m}/q) & \mbox{if $(q,m) = 1$} \vspace{2mm}\\ 0 & \mbox{otherwise}.\end{array}\right.\end{equation}
This immediately leads, via Lemma \ref{tizenketto}, to the following evaluation of $\mmhat(\theta)$.
\begin{lemmaenv}\label{lem29}
Suppose that $(a,q) = 1$ and that $\theta \in \mathfrak{M}_{a,q}$. then
\begin{equation}\label{29.4}\mmhat(\theta)  =  \left\{\begin{array}{ll} \displaystyle\frac{\mu(q)}{\phi(q)}e\left(-\frac{ab\overline{m}}{q}\right)\tau\left(\theta - \frac{a}{q}\right) + O\left((\log N)^{-A}\right) & \mbox{if $(m,q) = 1$}\vspace{2mm}\\ O\left((\log N)^{-A}\right) & \mbox{otherwise}.\end{array}\right.\end{equation}
\end{lemmaenv}
\noindent\textit{The minor arcs.} In this subsection we look at $\mmhat(\theta)$ and $\mmmhat(\theta)$ when $\theta$ is not close to a rational with small denominator.
\begin{lemmaenv}\label{lem477} Suppose that $a,q$ are positive integers with $(a,q) = 1$, and let $\theta$ be a real number such that $|\theta - a/q| \leq 1/q^2$. Then 
\begin{equation}\label{eq71} \mmhat(\theta)  \ll  (\log N)^{10}\left(q^{-1/2} + N^{-1/5} +  N^{-1/2}q^{1/2}\right).\end{equation} Thus if $\theta \in \mathfrak{m}$ then $\mmhat(\theta) = O((\log N)^{-A})$.
\end{lemmaenv}
\noindent\textit{Remarks.} This is a well-known estimate, at least when $b = m = 1$. The first (unconditional) results of this type were obtained by I.M. Vinogradov, and nowadays it is possible to give a rather clean argument thanks to the identity of Vaughan \cite{vaughan}. Chapter 24 of Davenport's book \cite{davenport} describes the use of Vaughan's identity in the more general context of the estimation of sums $\sum_{n \leq N} \Lambda(n) f(n)$. To obtain Lemma \ref{lem477} we used this approach, but could afford to obtain results which are rather non-uniform in $m$ due to the restriction $m \leq \log N$ under which we are operating. Details may be found in the supplementary document \cite{GreenMinorArcs}. We remark that existing results in the literature concerning minor arcs estimates for primes restricted to arithmetic progressions, such as \cite{BalogPerelli,Lavrik}, strive for a much better dependence on the parameter $m$.\endproof
\begin{lemmaenv}
Suppose that $a,q$ are positive integers with $(a,q) = 1$, and let $\theta$ be a real number such that $|\theta - a/q| \leq 1/q^2$. Then 
\begin{equation}\label{eq71b} \mmmhat(\theta)  \ll  (\log N)^3 \left(q^{-1} + qN^{-1}  + N^{-1/8A}\right).\end{equation} Thus if $\theta \in \mathfrak{m}$ then $\mmmhat(\theta) = O((\log N)^{-A})$.
\end{lemmaenv}
\proof Let $p_1,\dots,p_k$ be the primes less than or equal to $Q$ which do not divide $m$. Another application of the inclusion-exclusion principle gives
\[ \mmmhat(\theta)  =  N^{-1}e(-b\theta/m)\prod_{i = 1}^k\left(1 - \frac{1}{p_i}\right)^{-1}h(\theta), \]
where
\begin{equation}\label{eq028} h(\theta)  =  
\sum_{s = 0}^k (-1)^s \sum_{1 \leq i_1 < \dots < i_s \leq k} \;\; \sum_{\substack{1 \leq y \leq Nm/p_{i_1}\dots p_{i_s} \\ y \equiv b\mdsub{m}}} e\left(\frac{\theta p_{i_1}\dots p_{i_s} y}{m}\right).\end{equation}
Summing the geometric progression, one sees that the inner sum is no more than
\[ \min \left\{\Vert \theta p_{i_1}\dots p_{i_s}\Vert^{-1},2mN/p_{i_1}\dots p_{i_s}\right\}.\] 
We will split the sum over $s$ in \eqref{eq028} into two pieces, over the ranges $s \in [0,t]$ and $s \in (t,k]$ where $t = \log N/2A\log\log N$. Each of the primes $p_i$ is at most $Q \leq (\log N)^A$, so the product of any $s \leq t$ of them is no more than $\sqrt{N}$. Of course, all such products are distinct and so
\[ \sum_{s = 0}^t (-1)^s \sum_{1 \leq i_1 < \dots < i_s \leq k} \sum_{\substack{y \leq Nm/p_{i_1}\dots p_{i_s}\\y \equiv b\mdsub{m}}} e\left(\frac{\theta p_{i_1}\dots p_{i_s} y}{m}\right)  \leq  \sum_{n \leq \sqrt{N}} \min(\Vert \theta n \Vert^{-1},2mN/n).\]
This is a quantity whose estimation is standard in this area because of its pertinence to the estimation of exponential sums on minor arcs. It is bounded above by $C(\log N)^3(N^{1/2} + q + Nq^{-1})$; details may once again be found in \cite{GreenMinorArcs}.

On the other hand
\begin{eqnarray*}
& & \sum_{s = t+1}^k (-1)^s \sum_{1 \leq i_1 < \dots < i_s \leq k} \sum_{\substack{y \leq Nm/p_{i_1}\dots p_{i_s}\\ y \equiv b \mdsub{m}}} e\left(\frac{\theta p_{i_1}\dots p_{i_s} y}{m}\right) \\ & \leq & 2mN\sum_{s = t+1}^k \sum_{1 \leq i_1 < \dots < i_s \leq k} \prod_{j=1}^s p_{i_j}^{-1} \\ & \leq & 2mN\sum_{s = t+1}^k (s!)^{-1}\left(p_1^{-1} + \dots + p_k^{-1}\right)^s \\ & \leq & 4mN(2e\log\log\log N/t)^t \leq  mN^{1-1/4A} \leq  N^{1-1/8A}.\end{eqnarray*}
Since $\prod_{i=1}^k (1 - 1/p_i)^{-1} \ll \log N$, the claimed bound follows.\endproof\vspace{11pt}

\noindent\textit{Proof of Proposition \ref{mainprop}.} Suppose first of all that $\theta \in \mathfrak{M}_{a,q}$ for some $a,q$, and recall Lemmas \ref{lem27} and \ref{lem29}. If $q$ is $Q$-smooth then 
\[ \left|\mmhat(\theta) - \mmmhat(\theta)\right|  =  O(N(\log N)^{-A}).\]
If $q$ is not $Q$-smooth then $q > Q$ and so we get
\begin{eqnarray*} \left|\mmhat(\theta) - \mmmhat(\theta)\right| & \leq & |\mmhat(\theta)| + |\mmmhat(\theta)| \\ & \leq & 2/\phi(q) + O((\log N)^{-A})\\ & \leq & 4\log\log Q/Q + O((\log N)^{-A}),\end{eqnarray*} the last estimate being contained in \cite{HardyWright}, Chapter 17. Since we are assuming that $Q \leq (\log N)^A$ this expression is $O(\log\log Q/Q)$. If, on the other hand, $\theta \in \mathfrak{m}$ then we have
\begin{eqnarray*} \left|\mmhat(\theta) - \mmmhat(\theta)\right| & \leq & |\mmhat(\theta)| + |\mmmhat(\theta)| \\ & = & O((\log N)^{-A}) \\& = & O(Q^{-1}).\end{eqnarray*}
This at last completes the proof of Proposition \ref{mainprop}.\endproof
\section{Restriction and majorant estimates for primes}\label{sec5} In this section we prove Theorems \ref{thmHLM} and \ref{primerest}.

We have already seen, in \eqref{later2} and \eqref{later3}, how Proposition \ref{mainprop} implies an $L^1$--$L^{\infty}$ estimate for the operator $f \mapsto f \ast \psi_j$ of the form \eqref{eq699}. In fact, we have
\begin{equation}\label{hh376} \Vert f \ast \psi_j \Vert_{\infty}  \ll  \frac{\log j}{2^j}\Vert f \Vert_1\end{equation} for all $j = 1,\dots,K+1$. For each fixed $j = 1,\dots,K$, one can use the Riesz-Thorin interpolation theorem to interpolate between \eqref{later1} and \eqref{hh376}. This theorem, which is discussed in \cite[Chapter 7]{GreenRKP}, is better known to analytic number theorists as the type of convexity principle that underpins many basic estimates on $\zeta$ and $L$-functions. It gives
\begin{equation}\label{uis2} \Vert f \ast \psi_j \Vert_p  \ll  j^{2/p}(\log j)^{1 - 2/p}2^{-(1 - 2/p)j} N^{-2/p}\Vert f \Vert_{p'}.\end{equation}
For $j = K+1$ another interpolation, now between \eqref{later22} and \eqref{hh376}, instead gives
\[ \Vert f \ast \psi_{K+1} \Vert_p  \ll  (\log N)^{2/p}(\log K)^{1 - 2/p}2^{-(1 - 2/p)K}.\]
Recalling at this point the definition \eqref{Kdef} of $K$ we see that this implies
\[ \Vert f \ast \psi_{K+1} \Vert_p  \ll  (\log N)^{-1/p}N^{-2/p}.\]
Summing this together with \eqref{uis2} for $j = 1,\dots,K$ gives, because of the decomposition \eqref{dyadic},
\[ \Vert f \ast \mm\Vert_p  \leq  C(p)N^{-2/p}\Vert f\Vert_{p'}.\]
As we have already remarked, Theorem \ref{primerest} follows by the principle of $T$ and $T^{\ast}$.\endproof\vspace{11pt}

Now we prove Theorem \ref{thmHLM}. Although we will need a slightly different result later on, this theorem seems to be the most elegant way to state the majorant property for the primes.\vspace{11pt}

\noindent\textit{Proof of Theorem \ref{thmHLM}.} Let $(a_n)_{n \in \mathcal{P}_N}$ be any sequence of complex numbers with $|a_n| \leq 1$ for all $n$. We apply Theorem \ref{primerest} to the function $f$ defined by $f(n) = a_n/\log n$. Writing out the conclusion of Theorem \ref{primerest} gives, for any $p > 2$,
\[ \int \left|\sum_n f(n)\log n e(n\theta) \right|^p\,d\theta  \ll_p  N^{p/2-1}\left( \sum_n f(n)^2 \log n\right)^{p/2}.\]
Therefore
\begin{eqnarray*}
\int \left| \sum_{n \in \mathcal{P}_N} a_ne(n\theta) \right|^p \, d\theta & \ll_p & N^{p/2-1}\left(\sum_{n \in \mathcal{P}_N} \frac{|a_n|^2}{\log n}\right)^{p/2} \\ & \ll_p & N^{p-1}(\log N)^{-p}.\end{eqnarray*}
However it is an easy matter to check that
\[ \int \left| \sum_{n \in \mathcal{P}_N} e(n\theta) \right|^p \, d\theta  \geq  \int_{|\theta| \leq 1/2N} \left| \sum_{n \in \mathcal{P}_N} e(n\theta) \right|^p \, d\theta  \gg  N^{p-1}(\log N)^{-p}.\]
This proves Theorem \ref{thmHLM} for $p > 2$. For $p = 2$ it is trivial using Parseval's identity.\endproof
\section{Roth's theorem in the primes} Let $A_0$ be a subset of the primes with positive relative upper density. By this we mean that there is a positive constant $\alpha_0$ such that, for infinitely many integers $n$, we have \begin{equation}\label{eq3279}|A \cap \mathcal{P}_n|  \geq  \alpha_0 n/\log n.\end{equation} This is not a particularly convenient statement to work with, and our first lemma derives something more useful from it.
\begin{lemmaenv}\label{tech}
Suppose that there is a set $A_0 \subseteq \mathcal{P}$ with positive relative density, but which contains no \emph{3AP}s. Then there is a positive real number $\alpha$ and infinitely many primes $N$ for which the following is true. There is a set $A \subseteq \{1,\dots,\lfloor N/2\rfloor\}$, and an integer $W \in [\frac{1}{8}\log\log N,\frac{1}{4}\log\log N]$ such that
\begin{itemize}
\item $A$ contains no \emph{3AP}s
\item $\mm(A) \geq \alpha$ for some $b$ with $(b,m) = 1$, where $m = \prod_{p \leq W}p$.
\end{itemize}
\end{lemmaenv}
\proof Take any $n \geq \alpha_0^{-3}$ for which \eqref{eq3279} holds. Let $W = \lfloor \frac{1}{4}\log\log n\rfloor$, and set $m = \prod_{p \leq W}p$. Choose $N$ to be any prime in the range $(2n/m,4n/m]$. Now there are certainly no more than $m$ elements of $A_0$ which share a factor with $m$, and no more than $n^{3/4}$ elements $x \in A_0$ with $x \leq n^{3/4}$. Thus
\[ \sum_{b : (b,m) = 1}\sum_{\substack{x \leq n \\ x \equiv b\mdsub{m}}} A_0(x) \log x  \geq  \alpha_0 n/2,\]
and for some choice of $b$ we have
\begin{equation}\label{vantroost} \sum_{\substack{x \leq n \\ x \equiv b\mdsub{m}}} A_0(x) \log x  \geq  \alpha_0 n/2\phi(m).\end{equation}
Write $A = m^{-1}\left((A_0 \cap[n]) - b\right)$. This set, being a part of $A_0$ subjected to a linear transformation, contains no $3$-term AP. It is also clear that $A \subseteq \{1,\dots,\lfloor N/2\rfloor\}$. Furthermore \eqref{vantroost} is equivalent to
\[ \sum_{\substack{x \leq N \\ mx + b \, \mbox{\scriptsize is prime}}} A(x)\log (mx + b)  \geq  \alpha_0 n/2\phi(m),\]
which implies that $\mm(A) \geq \alpha_0 n/2mN \geq \alpha_0/8$. The lemma follows, with $\alpha = \alpha_0/8$.\endproof\vspace{11pt}
The reason we stipulate that $A$ be contained in $\{1,\dots,\lfloor N/2\rfloor\}$ is that $A$ does not contain any 3APs when considered \textit{as a subset of} $\mathbb{Z}_N = \mathbb{Z}/N\mathbb{Z}$. This allows us to make us of Fourier analysis on $\mathbb{Z}_N$. If $f : \mathbb{Z}_N \rightarrow \mathbb{C}$ is a function we will write, for any $r \in \mathbb{Z}_N$,
\[ \widetilde{f}(r)  =  \sum_{x \in \mathbb{Z}_N} f(x)e(-rx/N).\] Observe that $f$ may also be considered as a function on $\mathbb{Z}$ via the embedding $\mathbb{Z}_N \hookrightarrow [N]$, and then $\widetilde{f}(r) = f^{\wedge}(r/N)$.

For notational simplicity write $\mu = \lambda_{b,m,N}$. We will consider $A$ and $\mu$ as functions on $\mathbb{Z}_N$. Write $a = A\mu$. We will continue to abuse notation by using $\mu$ and $a$ as measures. Thus, for example, $a(\mathbb{Z}_N) \geq \alpha$.

Now if $A$ contains no (non-trivial) 3APs then 
\begin{eqnarray}\nonumber \sum_{x,d} a(x)a(x+d)a(x+2d) & = & \sum_x a(x)^3 \\ \nonumber & \leq & \sum_x \mu(x)^3 \\ \label{eq722}& \leq & (\log N)^3/N^2.\end{eqnarray}
We are going to show that this forces $\alpha$ to be small. We will do this by constructing a new measure $a_1$ on $\mathbb{Z}_N$ which is \textit{set-like}, which means that $a_1$ behaves a bit like $N^{-1}$ times the characteristic function of a set of size $\sim \alpha N$. The new measure $a_1$ will be fairly closely related to $a$, and in fact we will be able to show that 
\begin{equation}\label{3300} \sum_{x,d} a_1(x)a_1(x+d)a_1(x+2d) \qquad  \mbox{is small}.\end{equation} This, it turns out, is impossible; an argument of Varnavides based on Roth's theorem tells us that a dense subset of $\mathbb{Z}_N$ contains \textit{lots} of 3APs. We will adapt his argument in a trivial way to show that the same is true of set-like measures.

The arguments of this section, then, fall into two parts. First of all we must define $a_1$, define the notion of ``set-like'' and then show that $a_1$ is indeed set-like. The key ingredient here is Lemma \ref{lem63}, which says that $\widetilde{\mu}$ is small away from zero. Secondly, we must formulate and prove a result of the form \eqref{3300}. For this we need Theorem \ref{primerest}, the restriction theorem for primes.

The idea of constructing $a_1$, and the technique for constructing it, has its origins in the notions of \textit{granularization} as used in a paper of I.Z. Ruzsa and the author \cite{GreenRuzsa}. In the present context things look rather different however and, in the absence of anything which might be called a ``grain'', we think the terminology of \cite{GreenRuzsa} no longer appropriate.

Let us proceed to the definition of $a_1$. Let $\delta \in (0,1)$ be a real number to be chosen later, and set
\[ R  =  \left\{r \in \mathbb{Z}_N : |\widetilde{a}(r)| \geq \delta\right\}.\]
Let $k = |R|$, and write $R = \{r_1,\dots,r_k\}$. Let $\epsilon \in (0,1)$ be another real number to be chosen later, and write $B(R,\epsilon)$ for the Bohr neighbourhood
\[ \left\{x \in \mathbb{Z}_N : \left\Vert \frac{xr_i}{N} \right\Vert \leq \epsilon  \; \forall i \in [k]\right\}.\]
Write $B = B(R,\epsilon)$ and set $\beta(x) = B(x)/|B|$. Define
\begin{equation}\label{snk22} a_1  =  a \ast \beta \ast \beta.\end{equation}
It is easy to see that \begin{equation}\label{eq345}a_1(\mathbb{Z}_N)  \geq  \alpha.\end{equation} In Lemma \ref{lem117} below we will show that $\Vert a_1 \Vert_{\infty} \leq 2/N$, provided that a certain inequality between $\epsilon, k$ and $W$ is satisfied. This is what we mean by the statement that $a_1$ is set-like. 
\begin{lemmaenv}\label{lem63} Suppose that $N$, and hence $W$, is sufficiently large. We have
\[ \sup_{r \neq 0} |\widetilde{\mu}(r)|  \leq  2\log\log W/W.\]
\end{lemmaenv}
\proof Recall that $\widetilde{\mu}(r) = \mu^{\wedge}(r/N)$. There are three different cases to consider.\vspace{8pt}

\noindent Case 1. $r/N \in \mathfrak{M}_{0,1}$, that is to say $|r/N| \leq (\log N)^B/N$. Then by Lemma \ref{lem29} we have the asymptotic
\[ \widetilde{\mu}(r)  =  \tau(r/N) + O(\log N)^{-A}.\]
Observe, however, that $\tau(r/N) = 0$ provided that $r \neq 0$.\vspace{8pt}

\noindent Case 2. $r/N \in \mathfrak{M}_{a,q}$. Then Lemma \ref{lem29} gives
\[ \widetilde{\mu}(r)  =  \frac{\chi_q\mu(q)}{\phi(q)}e\left(-\frac{ab \overline{m}}{q}\right)\tau\left(\frac{r}{N} - \frac{a}{q}\right) + O(\log N)^{-A},\] where 
\[ \chi_q  =  \left\{\begin{array}{ll}1 & (q,m) = 1 \vspace{2mm}\\ 0 & \mbox{otherwise}.\end{array}\right. \] Since $m = \prod_{p \leq W} p$, we certainly have $\chi_q = 0$ for $q \leq W$. Thus indeed
\[ |\widetilde{\mu}(r)|  \leq  \sup_{n \geq W}\phi(n)^{-1} + O(\log N)^{-A}  \leq  2\log\log W/W.\] 
\noindent Case 3. $r/N \in \mathfrak{m}$. Then Lemma \ref{lem477} gives $\widetilde{\mu}(r) = \mu^{\wedge}(r/N) = O((\log N)^{-A})$.\endproof
\begin{lemmaenv}\label{lem117} Suppose that $\epsilon^k \geq 2\log\log W/W$. Then the measure $a_1$ is set-like, in the sense that we have $\Vert a_1 \Vert_{\infty} \leq 2/N$.
\end{lemmaenv}
\proof Indeed
\begin{eqnarray*}  a_1(x) & = & a \ast \beta \ast \beta(x) \\  & \leq & \mu \ast \beta \ast \beta(x) \\  & = & N^{-1}\sum_r \widetilde{\mu}(r)\widetilde{\beta}(r)^2e(rx/N) \\  & \leq &  N^{-1}\widetilde{\mu}(0)\widetilde{\beta}(0)^2 + N^{-1}\sum_{r \neq 0}|\widetilde{\mu}(r)||\widetilde{\beta}(r)|^2\\  & \leq & N^{-1} + N^{-1}\sup_{r \neq 0}|\widetilde{\mu}(r)|\sum_r |\widetilde{\beta}(r)|^2 \\ & = & N^{-1} + |B|^{-1}\sup_{r \neq 0}|\widetilde{\mu}(r)| \\ & \leq & N^{-1} + \frac{2\log\log W}{W|B|}.\end{eqnarray*}
Now by a well-known application of the pigeonhole principle we have $|B| \geq \epsilon^k N$, from which the lemma follows immediately.\endproof\vspace{11pt}

We move on now to the second part of our programme, which is a statement and proof of a result of the form \eqref{3300}.
\begin{propositionenv}\label{lem45} We have
\[ \sum_{x,d} a_1(x)a_1(x+d)a_1(x+2d)  \leq   C'N^{-3/2} + \frac{1}{N}\left(2^{12}\epsilon^2 \delta^{-5/2}+ C\delta^{1/2}\right).\]
\end{propositionenv} We will require several lemmas. The most important is a ``discrete majorant property''. Before we state and prove this, we give an elegant argument of Marcinkiewicz and Zygmund \cite{Zyg}. We outline the argument here since we like it and, possibly, it is not particularly well-known.
\begin{lemmaenv}[Marcinkiewicz--Zygmund]\label{MZlem} Let $N$ be a positive integer, and let $f : [N] \rightarrow \mathbb{C}$ be any function. Consider $f$ also as a function on $\mathbb{Z}_N$. Let $p > 1$ be a real number. Then we have
\[ \sum_{r \in \mathbb{Z}_N} |\widetilde{f}(r)|^p  =  \sum_{r = 0}^{N-1} |f^{\wedge}(r/N)|^p  \leq  C(p)N\int |\widetilde{f}(\theta)|^p\, d\theta.\]
\end{lemmaenv}
\proof Consider the function 
\[ g(n)  =  2\left(1 - \frac{|n|}{2N}\right)\chi_{|n| \leq 2N} - \left(1 - \frac{|n|}{N}\right)\chi_{|n| \leq N}.\] This function is equal to $1$ for all $n$ with $|n| \leq N$. Its Fourier transform, $g^{\wedge}(\theta)$, is equal to $2K_{2N}(\theta) - K_{N}(\theta)$, a difference of two Fej\'er kernels. Thus we have 
\[ f^{\wedge}  =  f^{\wedge} \ast \left(2K_{2N} - K_{N}\right),\]
and so
\begin{eqnarray*}
|\widetilde{f}(r)|^p & = & |f^{\wedge}(r/N)|^p \\& = & \left| \int f^{\wedge}(\theta)\left(2K_{2N}( r/N - \theta) - K_N( r/N - \theta)\right)\, d\theta \right|^p \\ & \leq & 3^{p-1}\left( 2^{p}\left|\int f^{\wedge}(\theta) K_{2N}( r/N - \theta)\,d\theta \right|^p + \left|\int f^{\wedge}(\theta) K_{N}( r/N - \theta)\,d\theta \right|^p\right) \\ & \leq & 3^{p-1}\left( 2^p\int |f^{\wedge}(\theta)|^p K_{2N}( r/N - \theta)\,d\theta  + \int |f^{\wedge}(\theta)|^p K_{N}( r/N - \theta)\,d\theta  \right)\end{eqnarray*}
by two applications of Jensen's inequality. It is necessary, of course, to use the fact that the Fej\'er kernels are non-negative. To conclude the proof, one only has to show that
\[ \sum_{r = 0}^{N-1} K_{N}( r/N - \theta)  \leq  CN,\]
together with a similar inequality for $K_{2N}$. But this is a straightforward matter using the bound
\[ \sum_{r = 0}^{N-1} K_{N}( r/N - \theta)  \leq  \sum_{j=0}^{N-1} \sup_{\phi \in [\frac{j}{N}, \frac{j+1}{N}]} K_{N}(\phi)\] together with the estimate
\[ K_N(\phi)  \ll  \min(N,N^{-1}|\phi|^{-2}),\]
valid for $|\phi| \leq 1/2$.\endproof
\begin{lemmaenv}[Discrete majorant property]\label{lem711} Suppose that $p > 2$. Then there is an absolute constant $C(p)$ \emph{(}not depending on $a$\emph{)} such that
\[ \sum_{r} |\widetilde{a}(r)|^p  \leq  C(p).\]
\end{lemmaenv}
\proof A direct application of Theorem \ref{primerest} gives
\[ \int |a^{\wedge}(\theta)|^p \, d\theta  \leq  C'(p)N^{-1}.\] The lemma is immediate from this and Lemma \ref{MZlem}.\endproof
\begin{lemmaenv}\label{lem8}
Suppose that $r \in R$. Then $\left|1 - \widetilde{\beta}(r)^4\widetilde{\beta}(-2r)^2\right| \leq 2^{12}\epsilon^2$.
\end{lemmaenv}
\proof We have 
\begin{eqnarray*}
\left|1 - \widetilde{\beta}(r)\right| & = & \frac{1}{|B|} \left|\sum_{x \in B} \left(1 - e(rx/N)\right)\right| \\ & = & \frac{1}{|B|} \left|\sum_{x \in B} \left(1 - \cos(2\pi r x/N)\right)\right| \\ & \leq & 4\pi^2 \sup_{x \in B} \Vert rx/N\Vert^2 \\ & \leq & 16\epsilon^2.\end{eqnarray*}
A very similar calculation shows that
\[ \left|1 - \widetilde{\beta}(-2r)\right|  \leq  64\epsilon^2,\] and the lemma follows quickly.\endproof\vspace{11pt}

\noindent\textit{Proof of Proposition \ref{lem45}.} By \eqref{eq722} we have, observing that $\widetilde{a_1} = \widetilde{a}\widetilde{\beta}^2$,
\begin{eqnarray}\nonumber
\sum a_1(x)a_1(x+d)a_1(x+2d) & \leq & \sum a_1(x)a_1(x+d)a_1(x+2d) \\ \nonumber & & \qquad - \sum a(x)a(x+d)a(x + 2d) + (\log N)^3N^{-2} \\ & = & O(N^{-3/2}) \nonumber\\ & & \; - N^{-1}\sum_r \widetilde{a}(r)^2 \widetilde{a}(-2r)\left(1 - \widetilde{\beta}(r)^4\widetilde{\beta}(-2r)^2\right).\label{eq654}\end{eqnarray}
Split the sum in \eqref{eq654} into two parts, that over $r \in R$ and that over $r \notin R$. When $r \in R$ we use Lemma \ref{lem8} to get
\begin{eqnarray*}
\sum_{r \in R}\widetilde{a}(r)^2 \widetilde{a}(-2r)\left(1 - \widetilde{\beta}(r)^4\widetilde{\beta}(-2r)^2\right) & \leq & 2^{12}\epsilon^2|R| \\ & \leq & C\epsilon^2\delta^{-5/2},\end{eqnarray*} this last inequality following from Lemma \ref{lem711} with $p = 5/2$. To estimate the sum over $r \notin R$, we again use Lemma \ref{lem711} with $p = 5/2$. Indeed using H\"older's inequality we have
\begin{eqnarray*}
\left|\sum_{r \notin R} \widetilde{a}(r)^2 \widetilde{a}(-2r)\left(1 - \widetilde{\beta}(r)^4\widetilde{\beta}(-2r)^2\right)\right| & \leq & 2\sup_{r \notin R} |\widetilde{a}(r)|^{1/2}\sum_r |\widetilde{a}(r)|^{5/2} \\ & \leq & C\delta^{1/2}.\end{eqnarray*}
This concludes the proof of Proposition \ref{lem45}.\endproof \vspace{11pt}

By \eqref{eq345} and Lemma \ref{lem117}, $a_1$ behaves a bit like a measure associated to a set of size $\alpha N$. As promised, we use this information together with an argument originally due to Varnavides \cite{Var} to get a \textit{lower} bound on $\sum a_1(x)a_1(x+d)a_1(x+2d)$.
\begin{lemmaenv}\label{lem3}
For some absolute constant $C_2$ we have
\[ \sum_{x,d \in \mathbb{Z}_N} a_1(x)a_1(x+d)a_1(x + 2d)  \geq  \exp\left(-C_2\alpha^{-2}\log(1/\alpha)\right)N^{-1}.\]
\end{lemmaenv}
\proof Let $A' = \{x \in \mathbb{Z}_N : a_1(x) \geq \alpha/2N\}$. By Lemma \ref{lem117} we have
\[ \alpha  \leq   \sum a_1(x) \leq \frac{2|A'|}{N} + \frac{\alpha}{2N}|A^{\prime c}|,\] which implies that $|A'| \geq \alpha N/4$. We will give a lower bound for $Z$, the number of 3APs in $A'$. It is clear that $\sum a_1(x)a_1(x+d)a_1(x + 2d)$ is at least $\alpha^3 Z/8N^3$.

Now by Bourgain's theorem\footnote{We could equally well use Roth's original theorem here, at the expense of making any bounds for the relative density in Theorem \ref{mainthm} even worse.} \cite{Bourgain} there is a constant $C_1$ such that if \[ M   \geq  \exp\left(C_1\alpha^{-2}\log(1/\alpha)\right)\] then any subset of $\{1,\dots,M\}$ of density at least $\alpha/8$ contains a 3AP with non-zero common difference. Now there are exactly $N(N-1)$ non-trivial arithmetic progressions of length $M$ in $\mathbb{Z}_N$, and $A'$ will have density at least $\alpha/8$ on many of them. To estimate exactly how many, fix a common difference $d \neq 0$, and let $I = \{0,d,2d,\dots,(M-1)d\}$. We have $\sum_x A' \ast I(x) \geq \alpha NM/4$, but $A' \ast I(x) \leq M$ for every $x$. Thus another simple averaging argument shows that $A' \ast I(x) \geq \alpha M/8$ for at least $\alpha N/8$ values of $x$.

In total, then, there are at least $\alpha N^2/8$ progressions of length $M$ on which $A'$ has density at least $\alpha/8$. Each of them contains a 3AP consisting of elements of $A'$. No 3AP thus counted can arise from more than $M^2$ progressions of length $M$. Thus we have two different ways of bounding $Z$, and putting them together gives
\[ Z  \geq  \alpha N^2/8M^2.\] The lemma follows.\endproof\vspace{11pt}

Combining this with Proposition \ref{lem45}, we get
\begin{equation}\label{eq8} C'N^{-1/2} + 2^{12}\epsilon^2 \delta^{-5/2} + C\delta^{1/2}  \geq  \exp\left(-C_2\alpha^{-2}\log(1/\alpha)\right).\end{equation}
There are constants $C_3,C_4$ so that if we choose \[ \delta  =  \exp\left(-C_3\alpha^{-2}\log(1/\alpha)\right)\] and \[\epsilon  =  \exp\left(-C_4\alpha^{-2}\log(1/\alpha)\right)\] then \eqref{eq8} cannot hold, and we will have derived a contradiction to the assumption that $A$ contains no 3APs. We are permitted to choose any values of $\epsilon$ and $\delta$ so that the condition of Lemma \ref{lem3} is satisfied. Recalling that $k \leq C\delta^{-5/2}$ (a consequence of Lemma \ref{lem711}) and that $W \geq \log \log N/8$, we see that \eqref{eq8} can indeed be contradicted provided that
\begin{equation}\label{badbound} \alpha  \geq  C\sqrt{\frac{\log_5 N}{\log_4 N}}.\end{equation} The subscripts indicate the number of iterated logarithms, not the base to which those logarithms are taken!

Let us remind the reader of what it is that we have contradicted. We assumed that there was a subset $A_0 \subseteq \mathcal{P}$ of positive relative upper density, containing no 3AP. The number $\alpha$ was related to the relative upper density of $A_0$, via the slightly technical reductions made in Lemma \ref{tech}. A bound of the form \eqref{badbound} also holds for $\alpha_0$. That is, any subset of $\mathcal{P}_n$ with cardinality at least $Cn(\log_5 n)^{1/2}/\log n(\log_4 n)^{1/2}$ contains a 3AP.

By far the most important reason for us getting such a poor bound was the need to prove Lemma \ref{lem63}, which says that by passing to a subprogression of common difference $m = \prod_{p \leq W}p$ one can make the primes look somewhat uniform. This is a rather crude trick but we have not been able to get around it. Even if we could, the resultant bounds would surely be many miles from the probable truth, which is that \textit{any} subset of $[N]$ of cardinality $N(\log N)^{-1000}$ contains 3APs.

Let us conclude by remarking that the methods of this section use rather little about the primes. In fact by the same argument one could establish a Roth-type theorem relative to any measure $\mu : \mathbb{Z}_N \rightarrow \mathbb{R}^{+}$ for which one had good control on $\sup_{r \neq 0} |\widetilde{\mu}(r)|$ together with bounds for $\Vert \widetilde{f} \Vert_{p}$, for some $p \in (2,3)$ and any $f$ satisfying $0 \leq f(x) \leq \mu(x)$ pointwise. In practise bounds of this latter type will come by restriction theory arguments of the type given in \S \ref{sec5}. A more general setting for our arguments, along the lines just described, is given in \cite{green-tao-1}.
\section{Acknowledgements} The author would like to thank Tim Gowers for his insights into Vinogradov's three-primes theorem, which played a substantial part in the development of this paper. He would also like to thank Imre Ruzsa for helpful conversations, Jean Bourgain for drawing his attention to the references \cite{B1,B2} and the students who attended the course \cite{GreenRKP} for their enthusiasm.

\end{document}